\documentclass[amstex,epsf]{article}
\usepackage[dvips]{graphicx,color}
\usepackage{amssymb}
\usepackage{amsmath}
\usepackage{caption}
\makeatletter
\newtheorem{theorem}{Theorem}[section]         
\newtheorem{corollary}{Corollary}[section]
\newtheorem{lemma}{Lemma}[section]

\begin{document}

\title{Presentation of finite subgroups of mapping class group of genus 2 surface by Dehn twists} 
\author{Gou Nakamura and Toshihiro Nakanishi}
\date{}

\maketitle

\begin{abstract}
In this note we give presentations of all finite subgroups of the mapping class group of a closed surface of genus 2 by the Humphries generators up to conjugacy.
\end{abstract}


\section{Introduction}

Let  ${\cal MCG}_{g}$ denote the mapping class group of a closed orientable surface of genus $g$. The Dehn-Lickorish theorem states that  ${\cal MCG}_{g}$ is generated by Dehn twists along finitely many simple closed curves.  S.~P.~Humphries found $2g+1$ Dehn twists which generate ${\cal MCG}_{g}$ for $g>1$. A finite presentation in the Humphries generators was given by Wajnryb \cite{waj1983}.    

The objective of this note is to give a generator system represented by products of Humphries generators $\{\omega_j\}_{j=1}^5$ 
for all finite subgroups, up to topological conjugacy,  of the mapping class group ${\cal MCG}_{2}$ of genus $2$. In general, it is not easy to find finite subgroups of a group when only one of its presentations is given. For ${\cal MCG}_{g}$, due to the Nielsen realization theorem (S.~Kerckhoff \cite{kerck1983}),  its finite subgroup is represented by a group of holomorphic automorphisms on a closed Riemann surface of genus $g$.  S.~A.~Broughton \cite{bou1991} made a complete list of all groups, up to topological conjugacy, which arise as groups of holomorphic automorphisms on some closed Riemann surfaces of genus $2$. It is natural to ask how those groups are embedded in $\mathcal{MCG}_2$. 

Our main theorem represents generators by products of $\omega_i$ ($i=1,..., 5$) for each group in Broughton's list.  For its statement we introduce some notation.
Let a finite group $G$ act on a Riemann surface $R$ as a group of holomorphic automorphisms. If the genus of the factor surface $R/G$ is $h$ and the covering map $\pi: R\to  R/G$ is branched over $n$ points $p_1$,..., $p_n$ with branching orders $m_j$,
then $(h; m_1,...,m_n)$ is the {\it type} of the orbifold $R/G$. In stead of $(h;m_1,...,m_n)$, we often write $(h; \nu_1^{r_1},..., \nu_p^{r_p})$ 
if $\nu_j$ appears $r_j$ times in $(m_1,...,m_n)$.
Let $\zeta_0$,  $\zeta_1$, ..., $\zeta_4$ be as in the following table. 
They have the orders indicated in the table. 
 \begin{center}
\begin{tabular}{|c|l|}\hline
order  &   \\ \hline 
2 & $\zeta_0= \omega_1\omega_2\omega_3\omega_4\omega_5^2\omega_4\omega_3\omega_2\omega_1$   \\
6 & $\zeta_1  =  \omega_1\omega_2\omega_3\omega_4\omega_5$    \\
6 & $\zeta_2  =  \omega_1\omega_2\omega_4^{-1}\omega_5^{-1}$ \\
8 & $\zeta_3  =  \omega_1^2\omega_2\omega_3\omega_4$ \\
10 & $\zeta_4  =  \omega_1\omega_2\omega_3\omega_4$  \\ \hline
\end{tabular}
\end{center}
The list below shows the group $G_{\ast}$ corresponding to {\rm (2.$\ast$)} in {\rm  \cite{bou1991}}, the order $|G_{\ast}|$ and the orbifold type.
\begin{theorem} \label{thm1} A non-trivial finite subgroup of ${\cal MCG}_{2}$ of a closed orientable surface of genus $2$ is conjugate with one of the groups in the following list. 
\end{theorem}  
\begin{itemize}
\item[(2.a)]  $G_a=\langle x=\zeta_0 : x^2=1\rangle\cong \mathbb{Z}_2$, $2$, \ $(0;2^6)$.

\item[(2.b)]  $G_b=\langle x=\zeta_1^3 : x^2=1\rangle\cong \mathbb{Z}_2$, $2$, \ $(1;2^2)$.

\item[(2.c)]  $G_c=\langle x=\zeta_1^2: x^3=1\rangle\cong \mathbb{Z}_3$, $3$, \ $(0;3^4)$.

\item[(2.e)] $G_e=\langle x=\zeta_3^2 : x^4=1\rangle\cong \mathbb{Z}_4$, $4$, \ $(0;2^2,4^2)$.

\item[(2.f)]  $G_f=\langle x=\zeta_0, y=\zeta_1^3: x^2=y^2=[x, y]=1\rangle\cong \mathbb{Z}_2\times \mathbb{Z}_2$,
 $4$, \ $(0;2^5)$.

\item[(2.h)]  $G_h=\langle x=\zeta_4^2: x^5=1\rangle\cong \mathbb{Z}_5$, $5$, \ $(0;5^3)$.

\item[(2.i)] $G_i=\langle x=\zeta_1: x^6=1\rangle\cong\mathbb{Z}_6$, $6$, \ $(0;3, 6^2)$.

\item[(2.k1)] $G_{k1}=\langle x=\zeta_2: x^6=1\rangle\cong\mathbb{Z}_6$, $6$, \ $(0;2^2,3^2)$.

\item[(2.k2)] $G_{k2}=\langle x=\zeta_1^3, y=(\omega_4\zeta_2\omega_4^{-1})^2: x^2=y^3=1, xyx^{-1}=y^{-1}\rangle\cong D_3$, 
$6$, \ $(0;2^2,3^2)$.

\item[(2.l)] $G_l=\langle x=\zeta_3: x^8=1\rangle\cong \mathbb{Z}_8$,  $8$, \ $(0;2,8,8)$.

\item[(2.m)] $G_m=\langle x=\zeta_2^{-1}\zeta_3^2\zeta_2, y=\zeta_3^2: x^4=y^4=1, x^2=y^2, xyx^{-1}=y^{-1}\rangle\cong \tilde{D}_2$,  $8$, \ $(0;4,4,4)$.

\item[(2.n)] $G_n=\langle x, y : x^2=y^4=1, xyx^{-1}=y^{-1}\rangle\cong D_4$,  where $x=\zeta_1^3$, $y=(\omega_4\omega_2^{-1}\omega_1^{-1})\zeta_3^2(\omega_4\omega_2^{-1}\omega_1^{-1})^{-1}$, $8$, \ $(0;2^3,4)$.

\item[(2.o)] $G_o=\langle x=\zeta_4: x^{10}=1\rangle\cong \mathbb{Z}_{10}$, $10$, \ $(0;2,5,10)$.

\item[(2.p)] $G_p=\langle x=\zeta_0, y=\zeta_1: x^2=y^6=[x, y]=1\rangle\cong \mathbb{Z}_{2}\times\mathbb{Z}_{6}$, $12$, \ $(0;2;6,6)$.

\item[(2.r)] $G_r=\langle x, y: x^4=y^3=1, xyx^{-1}=y^{-1}\rangle\cong D_{4,3,-1}$, \\
where $x=(\omega_3\omega_5^{-1}\omega_1^{-1})\zeta_3^6(\omega_3\omega_5^{-1}\omega_1^{-1})^{-1},
 y=\zeta_1^4$, $12$, \ $(0;3, 4^2)$.

\item[(2.s)] $G_s=\langle x=\zeta_1^3, y=\omega_4\zeta_2\omega_4^{-1}: x^2=y^6=1,xyx^{-1}=y^{-1}\rangle\cong D_6$, 
 $12$, \ $(0;2^3,3)$.

\item[(2.u)] $G_u=\langle x, y: x^2=y^8=1,xyx^{-1}=y^{3}\rangle\cong D_{2,8,3}$, 
where $x=\zeta_1^3, y=(\omega_4\omega_2^{-1}\omega_1^{-1})\zeta_3(\omega_4\omega_2^{-1}\omega_1^{-1})^{-1}$, 16, \ $(0;2,4,8)$.

\item[(2.w)] 
$G_w=\left\langle x,y,z,w : \begin{array}{l}x^2=y^2=z^2=w^3=[y, z]=[y,w]=[z,w]=1 \\ xyx^{-1}=y, xzx^{-1}=zy, xwx^{-1}=w^{-1}\end{array}\right\rangle$, \\
where $x=(\omega_1\omega_2\omega_1)\zeta_1^3(\omega_1\omega_2\omega_1)^{-1}$, $y=\zeta_0, z=\zeta_1^3, w=\zeta_1^4$,  \\
$G_w\cong \mathbb{Z}_2\ltimes(\mathbb{Z}_2\times\mathbb{Z}_2\times\mathbb{Z}_3)$, 
$24$, \ $(0;2,4,6)$.

\item[(2.x)] $G_x=\langle x=\zeta_2^4, y=\zeta_3^2 : x^3=y^4=1, xy^2=y^2x, (xy)^3=1\rangle\cong SL_2(3)$, $24$, \ $(0;3^2,4)$

\item[(2.aa)] $G_{aa}=\left\langle x, y, u: \begin{array}{l}x^3=y^4=(xy)^3=1, xy^2=y^2x, u^2=xyx^{-1}y^2 \\ uxu^{-1}=y^{-1}x^{-1}y, uyu^{-1}=x^{-1}yx\end{array}\right\rangle$\\
where $x=\zeta_2^4, y=\zeta_3^2$, $u=(\omega_1\zeta_1\omega_4\omega_2\omega_1^{-1})\zeta_3(\omega_1\zeta_1\omega_4\omega_2\omega_1^{-1})^{-1}$, \\
$G_{aa}\cong GL_2(3)$,
$48$, \ $(0;2,3,8)$
\end{itemize}
We thank Susumu Hirose and Shigeru Takamura for their encouragement and valuable comments.

\section{Facts about mapping class groups}
\subsection{The mapping class group and its generators} 
Our basic reference for mapping class groups is Farb-Margalit's book \cite{farb2012}, in particular, Sections 3, 4 and 7. The mapping class group ${\cal MCG}_{g}$ is generated by the Humphries generators or Dehn twists $\omega_0$, $\omega_1$,...., $\omega_{2g}$ along the loops depicted in Figure 1 (See \cite[Theorem 4.14]{farb2012}.) We consider an auxiliary Dehn twist $\omega_{2g+1}$ along the loop $c_{2g+1}$. Then the following relations hold:
\begin{gather}
\omega_i\omega_j=\omega_j\omega_i \quad \mbox{if $|i-j|\geq 2$, $1\le i, j\le 2g+1$} \label{eq11}\\
\omega_i\omega_{i+1}\omega_i=\omega_{i+1}\omega_{i}\omega_{i+1} \qquad (1\le i\le 2g)  \label{eq12}\\
(\omega_1\omega_2\cdots\omega_{2g+1})^{2g+2}=1 \quad \mbox{(a chain relation)}  \label{eq13}\\
\mbox{If $\zeta_0=\omega_1\omega_2\cdots\omega_{2g}\omega_{2g+1}\omega_{2g+1}\omega_{2g}\cdots\omega_2\omega_1$, then $\zeta_0^2=1$} 
\label{eq14}
\end{gather}
\begin{figure}[htbp]
 \begin{center}
  \begin{tabular}{c}   \label{figure1}
   \begin{minipage}{0.33\hsize}
    \begin{center}
    \includegraphics[width= 60mm]{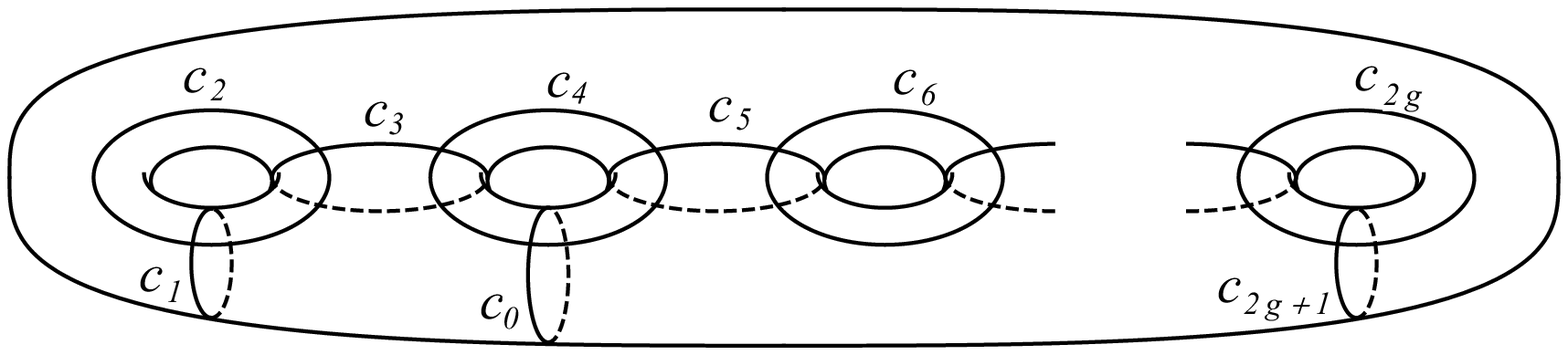}
    \caption{}
    \end{center}
   \end{minipage}\hspace{60pt} 
   \begin{minipage}{0.33\hsize}
    \begin{center}
    \includegraphics[width=35mm]{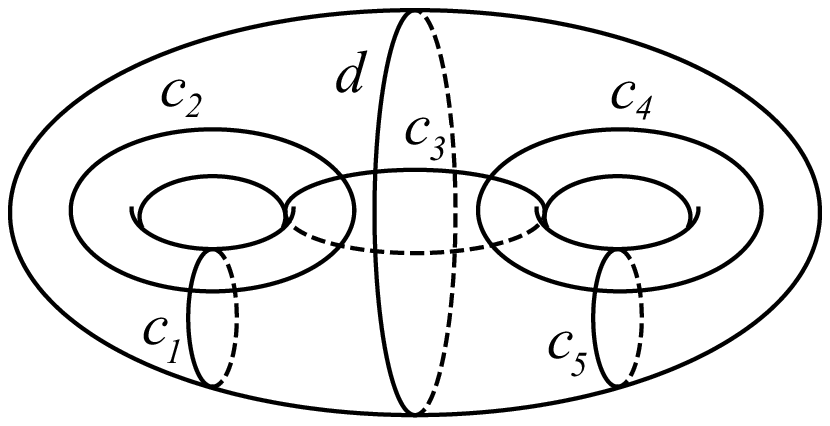}
    \caption{}
    \label{f2}
    \end{center}
   \end{minipage}
  \end{tabular}
 \end{center}
\end{figure}

\noindent
The element $\zeta_0$ is the {\it hyperelliptic involution}. If $1\le k\le m$, then  
\begin{equation}\label{eq20}
\begin{array}{l}
\omega_{i+k}(\omega_i\omega_{i+1}\cdots \omega_{i+m}) =(\omega_i\omega_{i+1}\cdots\omega_{i+m})\omega_{i+k-1}, \\
\\
\omega_{i+k}^{-1}(\omega_i^{-1}\omega_{i+1}^{-1}\cdots \omega_{i+m}^{-1}) =(\omega_i^{-1}\omega_{i+1}^{-1}\cdots\omega_{i+m}^{-1})\omega_{i+k-1}^{-1},
\end{array} 
\end{equation}
since we obtain from (\ref{eq11}) and (\ref{eq12}) 
\begin{eqnarray*}
\omega_{i+k}(\omega_i\omega_{i+1}\cdots \omega_{i+m}) &=& \omega_i\omega_{i+1}\cdots(\omega_{i+k}\omega_{i+k-1}\omega_{i+k})\omega_{i+k+1}\cdots\omega_{i+m}\\
&=&\omega_i\omega_{i+1}\cdots(\omega_{i+k-1}\omega_{i+k}\omega_{i+k-1})\omega_{i+k+1}\cdots\omega_{i+m}\\
&=&(\omega_i\omega_{i+1}\cdots\omega_{i+m})\omega_{i+k-1}.
\end{eqnarray*}
The second equation can be obtained in a similar way. Let
\[ \zeta=\omega_1\omega_2\cdots\omega_{2g+1}, \quad \eta=\omega_1\omega_2\cdots\omega_{2g}. \]
Applying (\ref{eq20})  we have $\omega_i\zeta=\zeta\omega_{i-1}$ for $i=2,...,2g+1$. Since 
\[\omega_1\zeta=\zeta\zeta^{-1}\omega_1\zeta=\zeta\omega_{2g+1}^{-1}\omega_{2g}^{-1}\cdots\omega_3^{-1}\omega_2^{-1}\zeta=\zeta^2\omega_{2g}^{-1}\cdots\omega_2^{-1}\omega_1^{-1},\]
we have $\omega_1\zeta=\zeta^2\eta^{-1}$,  $\omega_1\zeta^2=\zeta^2\omega_{2g+1}$ and $\omega_1\zeta^3=\zeta^2\omega_{2g+1}\zeta=\zeta^3\omega_{2g}$.
Continuing in this way we have for $i, j=1,..., 2g+1$,
\begin{equation}\label{eq21}
\omega_i\zeta^j=\zeta^j\omega_{i-j} \  (i\not=j), \quad  \omega_i\zeta^i=\zeta^{i+1}\eta^{-1},
\end{equation}
where the index $k$ for $\omega_k$ is considered modulo $2g+2$. From this follows
\begin{equation}\label{eq2a}
\omega_i=\zeta^{i+1}\eta^{-1}\zeta^{-i}, \quad (i=1,2,..., 2g+1).
\end{equation}

\begin{lemma}\label{lem22} {\rm (special cases of \cite[Proposition 4.12]{farb2012}) }
If $\eta=\omega_1\omega_2\cdots\omega_{2g}$ and $\xi=\omega_1^2\omega_2\cdots\omega_{2g}$, then $\eta^{2g+1}$ is a conjugate of $\zeta_0$ in {\rm (\ref{eq14})} and $\xi^{2g}=\eta^{2g+1}$.  Hence $\eta^{4g+2}=\xi^{4g}=1$ .
\end{lemma}

\noindent
{\it Proof.} Note that $\zeta_0=\zeta_0^{-1}$. From (\ref{eq2a}) 
\begin{eqnarray*}
\eta^{2g+1} &=& (\zeta^{-1}\omega_1^{-1}\zeta^2)(\zeta^{-2}\omega_2^{-1}\zeta^3)\cdots(\zeta^{-2g-1}\omega_{2g+1}^{-1})\\
 &=& \zeta^{-1}\omega_1^{-1}\omega_2^{-1}\cdots\omega_{2g+1}^{-1}\\
 &=& (\omega_{2g+1}\omega_{2g}\cdots\omega_2\omega_1\omega_1\omega_2\cdots\omega_{2g}\omega_{2g+1})^{-1}.
\end{eqnarray*}
A consequence of (\ref{eq20}) is $\omega_{i+1}\eta=\eta\omega_i$  ($i=1,..., 2g-1$). So we have  
\begin{equation*}
\xi^{2g} = \underbrace{\omega_1\eta\omega_1\eta\cdots \omega_1\eta}_{\mbox{$2g$ times}}
          = \omega_1\omega_2\cdots\omega_{2g}\eta^{2g}
          = \eta^{2g+1}.
\end{equation*}
From (\ref{eq2a})  we see that ${\cal MCG}_{g}$ is generated by three elements $\omega_0$, $\zeta$ and $\eta$. If $g=2$, then $\omega_0=\omega_5$. Therefore, 
we obtain a theorem by M.~Korkmaz for $g=2$.
\begin{corollary}{\rm (Korkmaz \cite{korkmaz2004})}
The mapping class group ${\cal MCG}_{2}$ is generated by $\zeta$ and $\eta$, where $\zeta^6=\eta^{10}=1$.
\end{corollary}

\subsection{Case of Genus $2$}

The mapping class group ${\cal MCG}_{2}$ is generated by Humphries generators $\omega_1$, $\omega_2$, $\omega_3$, $\omega_4$ and $\omega_5$ with defining relations (\ref{eq11}), (\ref{eq12}) and 
\begin{gather}
\zeta_1^6=(\omega_1\omega_2\omega_3\omega_4\omega_5)^6=1,  \nonumber\\
\zeta_0^2=(\omega_1\omega_2\omega_3\omega_4\omega_5^2\omega_4\omega_3\omega_2\omega_1)^2=1,  \nonumber\\
\mbox{$\omega_1\omega_2\omega_3\omega_4\omega_5^2\omega_4\omega_3\omega_2\omega_1$ and $\omega_i$ commute for $i=1,2,3,4,5.$}  \label{eq15}
\end{gather}
See \cite[p.184]{bir1974}.
Let $\zeta_0$,..., $\zeta_3$ and $\zeta_4$ be as in the table of Section 1.
From  (\ref{eq15}), the hyperelliptic involution $\zeta_0=\omega_1\omega_2\omega_3\omega_4\omega_5^2\omega_4\omega_3\omega_2\omega_1$ equals any conjugate of itself or its inverse. 
We proved in Lemma \ref{lem22} that $\zeta_3^4=\zeta_4^5=\zeta_0$, and hence $\zeta_3^8=\zeta_4^{10}=1$.  We shall show $\zeta_2^6=1$ in Lemma \ref{lem1245}, but topologically this arises from $2$-chain relations \cite[p.107]{farb2012} in one-holed tori on each side of the loop $d$ in Figure 2.
We remark that S.~Hirose studied in \cite{hirose2010} presentations of periodic mapping classes on orientable closed surfaces of genus $\le 4$ by Dehn twists by using topological and algebraic-geometric methods.

\section{Proof of the main theorem}
\subsection{Abelian Groups}  Our basic tools are elaborate  applications of (\ref{eq21}) for $g=2$ 
\begin{eqnarray}
\omega_i\zeta_1^j=\zeta_1^j\omega_{i-j} \  (i\not=j), \quad  \omega_i\zeta_1^i=\zeta_1^{i+1}\zeta_4^{-1},\label{eq30}
\end{eqnarray}
and trivial equations
\begin{eqnarray}
\zeta_1\omega_5=\zeta_0\omega_1^{-1}\omega_2^{-1}\omega_3^{-1}\omega_4^{-1}, \label{eq301}\\
\omega_1^{-1}\omega_2^{-1}\omega_3^{-1}\omega_4^{-1}\omega_5^{-1}=\zeta_1\zeta_0^{-1}=\zeta_1\zeta_0. \label{eq302}
\end{eqnarray}

\begin{lemma}\label{lem1245}
If $\zeta_2=\omega_1\omega_2\omega_4^{-1}\omega_5^{-1}$, then $\zeta_2^3=\zeta_0$. Hence $\zeta_2$ has order $6$.
\end{lemma}

\noindent
{\it Proof.}  The following equation deduced from (\ref{eq11}), (\ref{eq30}) and (\ref{eq302}) implies $\zeta_2^3=\zeta_0$.
\begin{eqnarray}
\zeta_2 &=& (\zeta_1\omega_5^{-1}\omega_4^{-1}\omega_3^{-1})(\omega_3\omega_2\omega_1\zeta_1\zeta_0)  =\zeta_1\omega_2\omega_5^{-1}\omega_1\omega_4^{-1}\zeta_1\zeta_0 \nonumber\\
&=& \omega_3(\zeta_1\omega_5^{-1})(\omega_1\zeta_1)\omega_3^{-1}\zeta_0 
= (\omega_3\zeta_4)\zeta_1^2(\omega_3\zeta_4)^{-1}\zeta_0. \label{eq1245}
\end{eqnarray}

Proofs of Theorem \ref{thm1} for cyclic groups  $G_c=\langle\zeta_1^2\rangle\cong \mathbb{Z}_3$, $G_e=\langle \zeta_3^2\rangle\cong \mathbb{Z}_4$,   
$G_h=\langle \zeta_4^2\rangle\cong \mathbb{Z}_5$,
$G_{l}=\langle \zeta_3\rangle\cong \mathbb{Z}_8$ and $G_{o}=\langle \zeta_4\rangle\cong \mathbb{Z}_{10}$ are straightforward. Since $\zeta_0$ is in the center of $\mathcal{MCG}_2$, proofs for abelian groups $G_f=\langle\zeta_0,\zeta_1^3\rangle\cong \mathbb{Z}_2\times \mathbb{Z}_2$ and $G_{p}=\langle \zeta_0, \zeta_1\rangle\cong\mathbb{Z}_2\times \mathbb{Z}_6$ are also easy. 

There are isomorphic pairs of groups $(G_a,G_b)$ and $(G_i, G_{k1})$. The types of orbifold for $G_a$ and $G_{k1}$ imply  that they must contain the hyperelliptic involution. So $G_a=\langle \zeta_0\rangle\cong \mathbb{Z}_2$, $G_b=\langle\zeta_1^3\rangle\cong \mathbb{Z}_2$ and $G_{i}=\langle \zeta_1\rangle\cong\mathbb{Z}_6$.The equation (\ref{eq1245}) means that $G_{k1}=\langle\zeta_2\rangle\cong\langle\zeta_1^2\rangle\times\langle\zeta_0\rangle\cong\mathbb{Z}_3\times\mathbb{Z}_2\cong\mathbb{Z}_6$. 
\subsection{Non-abelian Groups}

Groups with two or more generators require much effort. Assume that $x$ and $y$ generate a finite group of $\mathcal{MCG}_2$. We know each of them is a conjugate of a power of some $\zeta_i$, ($i=0,1,..., 4$), but finding suitable conjugates of $x$ and $y$ so that their product has also a finite order is rather laborious. 

\subsubsection{Groups $G_{k2}$ and $G_s$}

Let $a=\zeta_1^3$ and $b=\omega_1\omega_2\omega_5^{-1}\omega_4^{-1}=\omega_4\zeta_2\omega_4^{-1}$. Then $a^2=b^6=1$. It holds $(ba)^2=baba^{-1}=1$ because (\ref{eq30}) yields
$$ba=\omega_1\omega_2\omega_5^{-1}\omega_4^{-1}\zeta_1^3=\omega_1\omega_2\zeta_1^3\omega_2^{-1}\omega_1^{-1}.$$
So $a$ and $b$ generate $G_s\cong D_6$, and $a$ and $b^2$ generate $G_{k2}\cong D_3$. 


\subsubsection{Groups $G_n$, $G_u$ }
The group $G_u=\langle x, y: x^2=y^8=1, xyx^{-1}=y^3\rangle$ has the subgroup
$$G_n=\langle x, z=y^2: x^2=z^4=1, xzx^{-1}=z^{-1}\rangle.$$
If $x=\zeta_1^3$ and 
$y=\omega_1\omega_2\omega_4\omega_3\omega_2=(\omega_4\omega_2^{-1}\omega_1^{-1})\zeta_3(\omega_4\omega_2^{-1}\omega_1^{-1})^{-1}$, 
then $x^2=1$, $y^4=\zeta_0$ and $y^8=1$. By using (\ref{eq15}), (\ref{eq30}) and (\ref{eq301}), 
\begin{eqnarray*}
xyx^{-1}=\zeta_1^3(\omega_1\omega_2\omega_4\omega_3\omega_2)\zeta_1^{-3}&=& \omega_4\omega_5\zeta_1(\zeta_1\omega_5)(\zeta_1\omega_3\omega_2)\zeta_1^{-3}\\
&=& \omega_4\omega_5\zeta_1(\zeta_0\omega_1^{-1}\omega_2^{-1}\omega_3^{-1}\omega_4^{-1})(\omega_4\omega_3)\zeta_1^{-2}\\
&=& \omega_4\omega_5\omega_2^{-1}\omega_3^{-1}\zeta_1^{-1}\zeta_0\\
&=& (\omega_1\omega_2\omega_4\omega_3\omega_2)^{-1}\zeta_0=(\omega_1\omega_2\omega_4\omega_3\omega_2)^{3}=y^3.
\end{eqnarray*}

\subsubsection{Groups $G_r$, $G_w$}

Let $G_w=\mathbb{Z}_2\ltimes(\mathbb{Z}_2\times\mathbb{Z}_2\times\mathbb{Z}_3)$ have the  presentation as in (2.w). Since $[z,w]=1$, $u=zw$ generates $\mathbb{Z}_6\cong\mathbb{Z}_2\times\mathbb{Z}_3$ and satisfies $u^3=z$ and $u^4=w$. Therefore $G_w$ can be written as
\[\langle x, y, u : x^2=y^2=u^6=[x,y]=[y,u]=1, xux^{-1}=u^{-1}y\rangle.\]
This is an extension of the abelian group $G_p=\langle y, u: y^2=u^6=[y,u]=1\rangle$.
%
By letting $t=xu^3$ and $w=u^4$, we find the subgroup
\[ G_r=\langle t, w : t^4=w^3=1, twt^{-1}=w^{-1} \rangle\]
%
of $G_w$. Let $\zeta_5=\omega_1\omega_2\omega_1\omega_4^{-1}\omega_5^{-1}\omega_4^{-1}$. By using (\ref{eq11}), (\ref{eq12}), (\ref{eq30})  and (\ref{eq302}) we have
\begin{eqnarray*}
\zeta_5 &=& \zeta_1\omega_5^{-1}\omega_4^{-1}\omega_3^{-1}\omega_1\omega_5^{-1}\omega_3\omega_2\omega_1\zeta_1\zeta_0 = \zeta_1\zeta_5\zeta_1\zeta_0\\
&=& \zeta_1\omega_2\omega_1\omega_2\omega_4^{-1}\omega_5^{-1}\omega_4^{-1}\zeta_1\zeta_0=\omega_3(\zeta_1\omega_4^{-1}\omega_1\omega_5^{-1}\omega_2\zeta_1)\omega_3^{-1}\zeta_0\\
&=& \omega_3\omega_5^{-1}(\zeta_1\omega_5^{-1}\omega_1\zeta_1)\omega_1\omega_3^{-1}\zeta_0\\
&=& (\omega_3\omega_5^{-1}\omega_1^{-1})\zeta_3^2(\omega_3\omega_5^{-1}\omega_1^{-1})^{-1}\zeta_0
=(\omega_3\omega_5^{-1}\omega_1^{-1})\zeta_3^6(\omega_3\omega_5^{-1}\omega_1^{-1})^{-1}.
\end{eqnarray*}
Thus we obtain
\begin{equation}\label{eq35}
\zeta_5^2=\zeta_0, \quad \zeta_5^{-1}\zeta_1\zeta_5=\zeta_1^{-1}\zeta_0.
\end{equation}
and also
\begin{equation}\label{eq351}
(\zeta_5\zeta_1^k)^2=\zeta_0 \ \mbox{if $k$ is even}, \ (\zeta_5\zeta_1^k)^2=1 \ \mbox{if $k$ is odd}.
\end{equation}
Let $x=\zeta_5\zeta_1^3=(\omega_1\omega_2\omega_1)\zeta_1^3(\omega_1\omega_2\omega_1)^{-1}$, $y=\zeta_0$ and $u=\zeta_1$. Then $y^2=u^6=[x, y]=[y, u]=1$. The relations $x^2=1$ and $xux^{-1}=u^{-1}y$ are equivalent to $(\zeta_5\zeta_1^3)^2=1$ and $(\zeta_5\zeta_1^4)^2=\zeta_0$, which follow from (\ref{eq351}).


\subsubsection{Groups $G_m$, $G_x$, $G_{aa}$}\label{groupx}
$SL_2(3)=\langle x, y : x^3=y^4=1, xy^2=y^2x, (xy)^3=1\rangle$, where
\begin{equation*}\label{eq37a}
x=\left(\begin{array}{cc}
1 & 1\\
0 & 1\end{array}\right), \quad y=\left(\begin{array}{cc}
0 & 1\\ -1 & 0\end{array}\right) \quad \mbox{(the entries are in $\mathbb{Z}/(3\mathbb{Z}))$}.
\end{equation*}
$GL_2(3)$ is obtained by adding to $SL_2(3)$ the matrix
\[u=\left(\begin{array}{cc}
0 & -1\\
-1 & -1\end{array}\right),
\]
which satisfies $u^8=1$, $uy^2=y^2u$ and 
\begin{equation}\label{eq38}
u^2=xyx^{-1}y^2, \quad uxu^{-1}=y^{-1}x^{-1}y, \quad uyu^{-1}=x^{-1}yx.
\end{equation}
By letting $v=x^{-1}yx$, we find the subgroup $G_m$ of $SL_2(3)$ presented by  
$$G_m=\langle v, y: v^4=y^4=1, v^2=y^2, vyv^{-1}=y^{-1}\rangle.$$
Now, we represent $x$, $y$ and $u$ by $\omega_1$,..., $\omega_5$. By using (\ref{eq11})  and (\ref{eq30})
\begin{eqnarray*}
\zeta_2\zeta_3^{-2} &=& (\omega_1\omega_2\omega_4^{-1}\omega_5^{-1})(\omega_4^{-1}\omega_3^{-1}\omega_2^{-1}\omega_1^{-2})(\omega_5\zeta_1^{-1}\omega_1^{-1})\\
&=& \omega_1\omega_2(\omega_4^{-1}\zeta_1^{-1}\omega_5)(\omega_1^{-1}\zeta_1^{-1})\omega_1^{-1}\\
&=& (\omega_1\omega_2)\zeta_1^{-2}(\omega_1\omega_2)^{-1}.
\end{eqnarray*}
Hence $(\zeta_2\zeta_3^{-2})^3=1$. Let $x=\zeta_2\zeta_0=\zeta_2^4$ and $y=\zeta_3^2$. Then we have  $x^3=y^4=1$, $y^2=\zeta_0$, $(xy)^3=(\zeta_2\zeta_3^{-2})^3=1$ and $xy^2=y^2x$.  
Let
$$\begin{array}{l}
a=\omega_1^{-1}\zeta_2\omega_1=\omega_2\omega_4^{-1}\omega_5^{-1}\omega_1=\omega_2\omega_1\omega_4^{-1}\omega_5^{-1},\\
b=\omega_1^{-1}\zeta_3\omega_1=\omega_1\omega_2\omega_3\omega_4\omega_1=\zeta_4\omega_1=\omega_2\omega_1\omega_2\omega_3\omega_4,\\
c=\omega_2\omega_3\omega_5\omega_4\omega_3=(\zeta_1\omega_4\omega_2^{-1}\omega_1^{-1})\zeta_3(\zeta_1\omega_4\omega_2^{-1}\omega_1^{-1})^{-1}
\end{array}$$
and $u=\omega_1c\omega_1^{-1}$. Then $u^8=1$. %
Since $(x,y,u)=(\omega_1a\zeta_0\omega_1^{-1}, \omega_1b^2\omega_1^{-1},\omega_1c\omega_1^{-1})$, 
the relations $uxu^{-1}=y^{-1}x^{-1}y$ and $b^2ca(b^2c) ^{-1}=a^{-1}$ are equivalent. 
By using (\ref{eq30}) we have $\omega_1\zeta_4=\zeta_4\omega_4^{-1}\omega_3^{-1}\omega_2^{-1}\zeta_4=\zeta_4^2\omega_3^{-1}\omega_2^{-1}\omega_1^{-1}$, and then
\begin{eqnarray*}
b^2c =\zeta_4(\omega_1\zeta_4\omega_1\omega_2\omega_3)\omega_5\omega_4\omega_3
= \zeta_4^3\omega_5\omega_4\omega_3
=\zeta_1^3. 
\end{eqnarray*}
On the other hand, from  $\omega_2\omega_1\omega_4^{-1}\omega_5^{-1}\zeta_1^3=\omega_2\omega_1\zeta_1^3\omega_1^{-1}\omega_2^{-1}$
we have 
$$(\omega_2\omega_1\omega_4^{-1}\omega_5^{-1}\zeta_1^3)^2=1.$$
Since $\zeta_1^3=\zeta_1^{-3}$ this means 
$b^2ca(b^2c)^{-1} = a^{-1}$.  The relations $uyu^{-1}=x^{-1}yx$ and $(ac)b^2(ac)^{-1}=b^2$ are equivalent. The last relation easily follows from
\begin{eqnarray*}\label{eq33}
ac = \omega_2\omega_1\omega_4^{-1}\omega_5^{-1}\omega_2\omega_3\omega_5\omega_4\omega_3
= \omega_2\omega_1\omega_4^{-1}\omega_2\omega_3\omega_4\omega_3
= \omega_2\omega_1\omega_4^{-1}\omega_2\omega_4\omega_3\omega_4 
= b.
\end{eqnarray*}
Finally we show the first relation $u^2=xyx^{-1}y^2$ in (\ref{eq38}), which is equivalent to $c^2=ab^2a^{-1}b^4$. Since $b^4=b^{-4}=\zeta_0$, 
$ab^2a^{-1}b^4=ab^{-2}a^{-1}=(ab^{-1}a^{-1})^2$. We obtain $c^2=(ab^{-1}a^{-1})^2$ from
\begin{eqnarray*}
ab^{-1}a^{-1}\zeta_0 &=& \omega_2\omega_1\omega_4^{-1}\omega_5^{-1}(\omega_1^{-1}\zeta_4^{-1}\zeta_0)\omega_5\omega_4\omega_1^{-1}\omega_2^{-1}\\
&=&\omega_2\underline{\omega_1}\omega_4^{-1}(\underline{\omega_1^{-1}}\omega_5\omega_4\omega_3\underline{\omega_2\omega_1})\omega_5\omega_4\underline{\omega_1^{-1}\omega_2^{-1}}\\
&=&\omega_2\omega_4^{-1}\omega_5\omega_4\omega_3\omega_5\omega_4
=\omega_2\omega_4^{-1}\omega_4\omega_5\omega_4\omega_3\omega_4\\
&=&\omega_2\omega_5\omega_3\omega_4\omega_3=c.
\end{eqnarray*}

\bigskip

{\sc 
Gou Nakamura : Science Division, Center for General Education, 
Aichi Institute of Technology,
1247 Yachigusa, Yakusa,  Toyota,
470-0392, Japan}

{\it E-mail address} : {\tt gou@aitech.ac.jp}\bigskip

{\sc 
Toshihiro Nakanishi : Department of Mathematics, Shimane University, 
Matsue, 690-8504, Japan}  

{\it E-mail address} : {\tt tosihiro@riko.shimane-u.ac.jp}

\end{document}